\documentclass[12pt,amscd]{amsart}
\newtheorem{thm}{Theorem}[section]
\usepackage[all]{xy}
\usepackage{graphicx}
\usepackage{amsmath,amsxtra,amssymb,latexsym, amscd,amsthm}
\usepackage{indentfirst}
\usepackage[mathscr]{eucal}
\usepackage{hyperref}
\hypersetup{colorlinks=true,linkcolor=red,filecolor=magenta,urlcolor=green}
\usepackage{tikz}
\usepackage{capt-of}
\usepackage{caption}

\newtheorem{cor}[thm]{Corollary}
\newtheorem{lem}[thm]{Lemma}
\newtheorem{prop}[thm]{Proposition}
\theoremstyle{definition}
\newtheorem{defn}[thm]{Definition}
 
\newtheorem{rem}[thm]{Remark}

\DeclareMathOperator{\Z}{\mathbb {Z}}

\DeclareMathOperator{\depth}{depth}

\DeclareMathOperator{\ass}{Ass}

\DeclareMathOperator{\pd}{pd}

\DeclareMathOperator{\Tor}{Tor}

\def\x {\mathbf x}
\def\KK {\mathbb K}

\def\mi {\mathfrak m}

\begin{document}
	
	\title[Vertex covers, Betti numbers and projective dimensions of binary trees]{The vertex covers, Betti numbers and projective dimensions of perfect binary trees}

	\author[Nguyen Thu Hang]{Nguyen Thu Hang }
	%\address{International Centre of Research and Postgraduate Training in %Mathematics, 18B Hoang Quoc Viet Street, Ha Noi, Vietnam}
	\address{Thai Nguyen University of Sciences, Phan Dinh Phung Ward, Thai Nguyen, Vietnam}
	\email{hangnt@tnus.edu.vn}
	
	\author[Tran Duc Dung]{Tran Duc Dung}
	%\address{International Centre of Research and Postgraduate Training in %Mathematics, 18B Hoang Quoc Viet Street, Ha Noi, Vietnam}
	\address{Thai Nguyen University of Sciences, Phan Dinh Phung Ward, Thai Nguyen, Vietnam}
	\email{dungtd@tnus.edu.vn}
	
	\author[Do Van Kien]{Do Van Kien}
	\address{Hanoi Pedagogical University 2, 32 Nguyen Van Linh Street, Xuan Hoa Ward, Phu Tho, Vietnam}
	\email{dovankien@hpu2.edu.vn}
	
	%\date{}

	%\commby{}
	%-----------------------------------------------------------
	\begin{abstract} Let $T$ be a perfect binary tree and $I$ be its edge ideal in the polynomial ring $S$. We determine the vertex cover number, independent number, and establish the recursive formula to compute the number of minimal vertex covers. As a consequence, we compute the depth and projective dimension of $S/I$ and show that the total Betti number of $S/I$ at the highest homological degree always equals one.
	\end{abstract}
	
	% -----------------------------------------------------------
	\maketitle
	% -----------------------------------------------------------
	\section{Introduction}
Throughout this paper, let $G$ be a finite simple graph with vertex set $V(G)=\{x_1,\dots, x_n\}$ and edge set $E(G)\subseteq V(G)\times V(G)$. Let $S = \mathbb K[x_1,\ldots,x_n]$ be a polynomial ring over a field $\mathbb K$. The edge ideal of $G$ is denoted by $I(G)=(x_ix_j |\ \{x_i,x_j\}\in E(G))\subseteq S$. A central theme in the study of edge ideals is understanding the relationships between the algebraic invariants of $I(G)$ and the combinatorial properties of the underlying graph $G.$ It shows many beautiful connections between algebra and combinatorics.  
A {\it vertex cover} of $G$ is a subset $C$ of $V(G)$ which meets every edge of $G$. A vertex cover is said to be {\it minimal} (with respect to set inclusion) if none of its proper subsets is itself a cover.  The {\it vertex cover number} of $G$ is defined as the smallest cardinality of a vertex cover set of $G$ and denoted by $\alpha(G)$. The number of minimal vertex covers of $G$ is denoted by $m(G)$.  A subset $W$ of $V(G)$ is called an {\it independent set} if no two vertices in $W$ are adjacent. The {\it independence number} is defined as the cardinality of a largest independent set of $G$ are denoted by $\beta(G)$. Obviously, $\alpha(G)=n-\beta(G)$. An independent set of $G$ is said to be {\it maximal} if it cannot be extended to a larger independent set. It is clear that a subset $W$ of $V(G)$ is a maximal independent set of $G$ if and only if $V(G)\setminus W$ is a minimal vertex cover of $G$. Therefore, $m(G)$ is also exactly the number of maximal independent sets of $G$.

Inspired by Erd\"{o}s and Moser who posed the problem of determining the maximum possible value of $m(G)$ in terms of the order of $G$, many authors have concentrated on studying $m(G)$ for various classes of graphs. We refer to \cite{Fu} for connected graphs, \cite{HT} for triangle-free graphs, \cite{Co84,Sa88,Wi86} for trees, \cite{Liu} for bipartite graphs, $\ldots$ Recently, Hoang and Trung in \cite[Theorem 2.7]{HT19} showed that $m(G)\le2^{\alpha (G)}$. Generally, obtaining an explicit formula for $m(G)$ is complicated. In \cite{B}, Bouchat gave an inductive formula for the number of minimal vertex covers of a path. However, for an arbitrary tree, the algorithms obtained are all quite complex time-based algorithms. In this paper, we investigate perfect binary trees. By combining combinatorial decomposition techniques with recursion on the height of the perfect binary tree, we derive a recurrence relation for computing this invariant.
\begin{thm}[Theorem \ref{numbercover}]
Let $T$ be a perfect binary tree of height $h$, and $m_h$ is the number of minimal vertex covers of $T$. Then there is a recursive relation
$$\begin{cases}
	m_0=1,m_1=2,m_2=4,\\
	m_{h+1}=2m_hm_{h-2}^4+m_{h-1}^4-m_{h-2}^8 \text{ for all } h\ge 2.
\end{cases}$$
\end{thm}

Minimal vertex covers encode primary decompositions and are closely related to the irredundant irreducible decompositions of the ideals (\cite{HSM}). It thus allows us to easily calculate the number of associated prime ideals of edge ideals of a perfect binary tree.
	
We also establish an explicit formula for the vertex cover number $\alpha(T)$.
	\begin{prop}[Proposition \ref{matnumber}] Let $T$ be a perfect binary tree of height $h$. Then 
		\begin{equation*}
			\alpha(T)=\dfrac{2^{h+2}-(3+(-1)^h)}{6}.
		\end{equation*}
	\end{prop}
\noindent This formula provides an effective way to compute the dimension of $S/I(T)$. It is natural to be concerned with the depth of $S/I(T)$, which is closely related to its projective dimension by Auslander-Buchsbaum's formula.

	By the definition, the projective dimension of $S/I(G)$ is the highest homological degree at which the Betti number is not zero. When $G$ is a tree, the homological structure of edge ideals is much better understood. Since every tree is a chordal graph, its edge ideal has a linear resolution by Fröberg’s Theorem, which imposes a particularly nice structure on the Betti table. In the $\Z^n$-graded setting, R. Bouchat \cite{B} proved that the nonzero finely graded Betti numbers of the edge ideal of a tree must be $1$. 
	
 It is not easy to describe the complete minimal resolution or Betti table of the edge ideal of a graph. For a perfect binary tree $T$, we rely on the result of calculating the depth of the quotient ring $S/I(T)$ via the vertex cover number of the tree to obtain the projection dimension of the ring $S/I(T)$. Consequently, we get the result of the projective dimension of $S/I(T)$ as follows:
	\begin{thm}[Theorem \ref{pdim}]
		Let $T$ be a perfect binary tree of height $h>0$ and $I$ be the edge ideal of $T$. Then the following assertions hold.
		\begin{enumerate}
			\item[(1)] If $h=3m$ then $\pd_S(S/I)=\dfrac{10(8^m-1)}{7}$.
			\item[(2)] If $h=3m+1$ then $\pd_S(S/I)=\dfrac{20.8^m-6}{7}$.
			\item[(3)] If $h=3m+2$ then $\pd_S(S/I)=\dfrac{40.8^m-5}{7}$.
		\end{enumerate}
	\end{thm}
By decomposing the perfect binary tree into subtrees and employing the long exact sequences arising from short exact sequences, we establish the following.
	\begin{thm}[Theorem \ref{main2}] Let $T$ be a perfect binary tree of height $h>0$ and $I$ its edge ideal. Denote by $p_h=\pd_S(S/I)$ the projective dimension of $S/I$. Then the last total Betti number of $S/I$ is $\beta_{p_h}(S/I) = 1.$
	\end{thm}
\noindent It therefore provides a class of non-Gorenstein rings of type 1. 
	
	\rm The paper is organized as follows. Section $2$ contains all the concepts that will be used throughout the paper. We distinguish here between combinatorial concepts arising from graph theory and notions from commutative algebra, and we describe the connections between these two areas.  In Section $3$, with a perfect binary tree $T$, we investigate the induction formula for finding the number of minimal covers of $T$. In particular, we found a specific formula for the vertex cover number of this graph. In Section $4$, we survey the minimal resolution of the edge ideal of the perfect binary tree. Moreover, by the calculation of the vertex cover number in Section $3$, we can deduce the depth of the quotient ring $S/I(T)$. From this, we can deduce the final Betti number of the free solution of the ideal $I(T)$.
	%%%%%%%%%%%%%%%%%%%%%%%%%%%%%%%%%%%%%%%%%%%%%%%%%%%%%%%%%%
	\section{Preliminaries}
	
\rm In this section, we recall some definitions and properties concerning minimal vertex covers of graphs, the Betti numbers, and the depth of edge ideals of graphs. The interested readers are referred to \cite{BH,D} for more details. Note that by abuse of notation, $x_i$ will at times be used to denote both a vertex of a graph and the corresponding variable of the polynomial ring.
	%Throughout the paper, let $\mathbb K$ be a field, and $R = \mathbb K[x_1,\ldots,x_n]$, with $n\geqslant 2$ be a polynomial ring, and let $\mi = (x_1,\ldots,x_n)$ be the maximal homogeneous ideal of $R$. 
	
	\subsection{Graphs and their edge ideals}
	
	Let $G$ be a simple graph, that is, it has no vertex connected to itself by an edge. We use the symbols $V(G)=\{x_1,\ldots, x_n\}$ and $E(G)=\{\{x_i,x_j\} | \ x_i, x_j\in  V(G)\}$ to denote the vertex set and the edge set of $G$, respectively. 
	
	We now recall several classes of graphs that we study in this work. Two vertices $x_i, x_j \in V(G)$ are called {\it adjacent} (or {\it neighbors}) if $\{x_i,x_j\} \in E(G)$. For a vertex $x_i$ of  $G$, we denote by $N(x_i)$ the set of all the neighbors of $x_i$, also called the {\it neighborhood} of  $x_i$. More precisely, $N(x_i) = \{x_j \in  V(G) :  \{x_i,x_j\} \in E(G)\}$. Moreover, let $N[x_i] = N(x)\cup \{x_i\}$ be the {\it closed neighborhood} of $x_i.$  A vertex $x_i$ is a {\it leaf} if $N(x_i)$ has cardinality one. A graph $H$ is called a {\it subgraph} of $G$ if $V(H)\subseteq V(G)$ and $E(H) \subseteq E(G)$.  A graph $H$ is called an {\it induced subgraph} of $G$ if the vertices of $H$ are vertices in $G$, and for vertices $x_i$ and $x_j$ in $V(H)$, $\{x_i,x_j\}$ is an edge of $H$ if and only if $\{x_i,x_j\}$ is an edge in $G$. For a subset $U$ of $V(G)$, we denote by $G\setminus U$  the induced subgraph of $G$ on $V(G)\setminus U$.
	
	A  graph is {\it connected} if there is a path from any vertex to any other vertex in the graph. A graph that is not connected is said to be {\it disconnected}. A connected component of a graph $G$ is a connected subgraph that is not part of any larger connected subgraph. 
	
	A {\it cycle} is a graph with $V(C_n)=\{x_1,\ldots,x_n\}$ and $V(G)=\{\{x_i,x_{i+1}\}| i=1,\ldots,n-1\}\cup \{x_n,x_1\}.$
	
	The graph $G$ is {\it bipartite} if  $V(G)$ can be partitioned into two subsets $X$ and $Y$ such that every edge has one end in $X$ and another end in $Y$, such a partition $(X, Y)$ is called a {\it bipartition} of the graph.  Note that $G$ is bipartite if and only if it has no cycle of odd length (see \cite[Theorem 4.7]{BM}).  A connected graph without cycles is a {\it tree}. Obviously, a tree is bipartite. In this survey, when the graph is a tree, we will always denote the graph by the symbol $T$. It is clear that a tree with $n$ vertices has $n-1$ edges.
	
The tree $T$ is called a {\it perfect binary tree} if all the leaf vertices are same height (the height of a vertex is defined as the number of edges from the root to a vertex), and all non-leaf vertices have two children. In simple terms, this means that all leaf vertices are at the maximum height of the tree. We denote the height of $T$ by $h$.
	\begin{center}
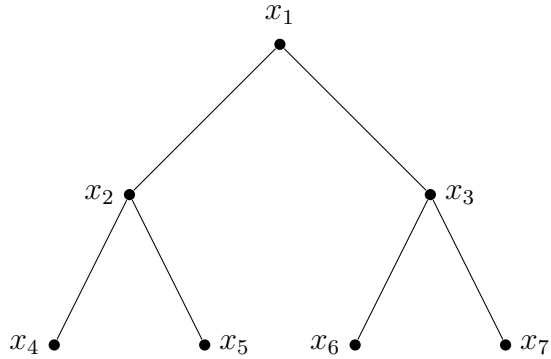

		\begin{tikzpicture}[every node/.style={circle, fill, inner sep=1.5pt},
			label distance=0mm]
			% Nodes
			\node (x1) at (0, 0) [label=above:$x_1$] {};
			\node (x2) at (-2, -2) [label=left:$x_2$] {};
			\node (x3) at (2, -2) [label=right:$x_3$] {};
			\node (x4) at (-3, -4) [label=left:$x_4$] {};
			\node (x5) at (-1, -4) [label=right:$x_5$] {};
			\node (x6) at (1, -4) [label=left:$x_6$] {};
			\node (x7) at (3, -4) [label=right:$x_7$] {};
			
			% Edges
			\draw (x1) -- (x2);
			\draw (x1) -- (x3);
			\draw (x2) -- (x4);
			\draw (x2) -- (x5);
			\draw (x3) -- (x6);
			\draw (x3) -- (x7);
			%\draw[dashed] (-3.2, -4.5) -- (3.2, -4.5);
		\end{tikzpicture}
		\vspace*{-0,4cm}
		\captionof{figure}{The perfect binary tree with height $h=2$}
	\end{center}
	
	\begin{lem}\label{nodes} Let $T$ be a perfect binary tree of height $h$. Then the number of leaves,  the number of vertices, and the number of edges of $T$ are $2^h$, $2^{h+1}-1$, and $2^{h+1}-2$, respectively.
	\end{lem}
	\begin{proof} The case where $h=0$ is trivial. We consider $h>0$ and assume that every perfect binary tree of height $h-1$ has $2^{h-1}$ leaves. Then by removing all the leaves of $T_h$ and the edges connected to them, we get a perfect binary subtree of height $h-1$. This subtree has $2^{h-1}$ leaves which implies that there are $2^{h-1}$ vertices of $T_h$ at level $h-1$. Each vertex has two children. Therefore, the number of leaves of $T_h$ is $2.2^{h-1}=2^h$. Hence, the number of vertices of $T$ is $1+2+\cdots + 2^h=2^{h+1}-1,$ and  the number of edges of $T$ is $(2^{h+1}-1)-1=2^{h+1}-2.$ 
	\end{proof}
\begin{defn}
Let $G$ be a finite simple graph with the vertex set $V(G)=\{x_1,\ldots,x_n\}$. The {\it edge ideal} of $G$ is defined by
	\begin{equation*}
		I(G) = (x_ix_j \mid \{x_i, x_j\} \in E(G)) \subseteq S=\mathbb K[x_1,\ldots,x_n].
	\end{equation*}
\end{defn}
Because $I(G)$ is a squarefree monomial ideal, minimal primes of $I(G)$ are variable ideals given by vertex covers. This means that the edge ideal $I(G)$ has the primary decomposition
\begin{equation} \label{intersect}
	I(G) = \bigcap \big(\x_{C}=\underset{x_i\in C}{\prod} x_i \mid C \text{ is a minimal vertex cover of } G\big).
\end{equation}
\noindent Observe that the height of $I(G)$ also is the minimal cardinality of a vertex cover of $G,$ that is $\alpha(G)$. It therefore implies the dimension of $S/I(G)$. 
	\begin{lem}[\cite{HM}, Lemma 1] $\dim S/I(G)=n-\alpha(G)=\beta (G).$
	\end{lem}
A subset $M=\{e_{1},\ldots, e_{s}\}\subseteq E(G)$ is called  {\it a matching} of $G$ if for all $i\neq j$, one has $e_i\cap e_j=\emptyset$. If $M$ is a matching, the two ends of each edge of $M$ are said to be matched under $M$, and each vertex incident with an edge of $M$ is said to be covered by $M$.  The number of edges in a maximum matching in a graph $G$ is called the {\it matching number} of $G$. 

In bipartite graphs, the matching numbers are closely related to vertex cover numbers by K\"{o}nig's theorem. 
\begin{lem}[\cite{D}]\label{konig} In any bipartite graph, the number of edges in a maximum matching equals the number of vertices in a minimum vertex cover.
\end{lem}
	\subsection{Depth and projective dimension} Let $S = \mathbb K[x_1,\ldots,x_n]$ be a polynomial ring over a field $\mathbb K$. The depth of graded modules and ideals over $S$ is one of the objects in our work. This invariant can be defined via either the minimal free resolutions or the local cohomology modules. 
	
	Let $M$ be a nonzero finitely generated graded  $S$-module.  Hilbert’s Syzygy Theorem confirms that $M$  has a minimal free resolution of length at most $n$.  Let
	$$0 \rightarrow \bigoplus_{j\in\Z} S(-j)^{\beta_{p,j}(M)} \rightarrow \cdots \rightarrow \bigoplus_{j\in\Z}S(-j)^{\beta_{0,j}(M)}\rightarrow 0$$
	be the minimal free resolution of $M$. The {\it projective dimension} of $M$ is defined as the length of this resolution and denoted by $\pd_S(M)$.
	The depth and the projective dimension of $M$  are related by Auslander-Buchsbaum's formula $$\depth(M)+\pd_S(M)=n.$$

	Note that the depth of $M$ can also be computed via the local cohomology modules of $M$. Let $H_{\mi}^i(M)$ be the $i$th cohomology module of $M$ with support in $\mi=(x_1,\ldots,x_n)$. Then,
	$$\depth(M): = \min\{i\mid H_{\mi}^i(M) \ne 0\}.$$
	
	%We have the following estimates on depth along a short exact sequence (see \cite[Proposition 1.2.9]{BH}).
	
	%\begin{lem}\label{exseq} Let $0 \rightarrow L \rightarrow  M \rightarrow N \rightarrow 0$ be a short exact sequence of finitely generated graded R-modules. Then
	%\begin{enumerate}
	%\item[(1)] $\depth(M)\geq \min\{\depth(L),\depth(N)\}.$
	%\item[(2)] $\depth(L)\geq \min\{\depth(M),\depth(N)+1\}.$
	%\end{enumerate}
	%\end{lem}
	%We recall the following well-known results. The first one is known in \cite[Corollary 1.3]{R}, and the second one is in \cite[Corollary 3.3]{R}.
	
	%\begin{lem}\label{lem_upperbound} Let $I$ be a monomial ideal and $f$ a monomial such that $f \notin I$. Then
	%\begin{enumerate}
	%\item[(1)] $\depth (S/I) \le \depth (S/(I:f)).$
	%\item[(2)] If $f$ is regular on $S/I$, then
	% $$\depth (S/I) = 1+\depth (S/(I,f)).$$
	%\end{enumerate}
	%  \end{lem}
%The following basic fact will be useful in several proofs. For clarity it is stated here.
%\begin{lem}{\rm (see \cite[Lemma 2.8]{SM})}\label{basicdepth}
%Let $I$ be an ideal in a polynomial ring $R$, let $x$ be an indeterminate over $R$, and let $S=R[x]$. Then $\depth S/IS=\depth R/I+1.$
%\end{lem}
\subsection{Minimal resolution and Betti numbers of edge ideals}
Let  $G = (V, E)$ be  a finite simple graph on a vertex set $V = \{x_1, x_2, \ldots, x_n\}$ and $I(G)$ the edge ideal of $G$. Since $I(G)$ is minimally generated by monomials of degree 2, the graded module $S/I(G)$ over the polynomial ring $S= \mathbb K[x_1, x_2, \ldots, x_n]$ has a unique minimal graded free resolution of length $p \leq n$ as follows:
\begin{align*}
	0 \to \bigoplus_{j=p+1}^{s_p} S[-j]^{\beta_{p,j}}
	\to \cdots \to
	\bigoplus_{j=i+1}^{s_i} S[-j]^{\beta_{i,j}}
	\to \cdots \to
	\bigoplus_{j=2}^{s_1} S[-j]^{\beta_{1,j}}&\to S \to\\
	&\to S/I(G) \to 0.
\end{align*}
By definition, the number $p=\pd_S(S/I(G))$ and the number $\beta_{i,j}$ of generators of $i$-th syzygy module in degree $j$ is called {\it the $i$-th graded Betti number} of $S/I(G)$ in degree $j$, denoted by $\beta_{i,j}(S/I(G))$ (or simply $\beta_{i,j}(G)$). It is well known that $\beta_{i,j}(S/I(G))=\dim_{\KK}\Tor_i^S(S/I(G),\KK)_j$. We also call $\beta_{i}(S/I(G))=\sum_j \beta_{i,j}(S/I(G))$ the $i$-th (total) Betti number of $S/I(G)$. 

If we consider $S$ as a finely graded ring (or $\Z^n$-graded ring), that is, $\deg(x_1^{a_1}\ldots x_n^{a_n})=(a_1,\ldots,a_n)=\mathbf a$, then the nonzero graded Betti numbers of the edge ideal of a tree must be $1$.
\begin{thm}[\cite{B}, Theorem 2.2.2]
	Let $T$ be a tree with $n$ edges, then the finely graded Betti numbers $\beta_{i,{\mathbf a}}(S/I(T))$ is either 0 or 1. 
\end{thm}
\section{Minimal vertex covers of perfect binary trees}
Let $G$ be a simple graph. It is well-known that the minimal primary decomposition of the edge ideal $I(G)$ of $G$ deeply relates to the set of all minimal vertex covers of $G$. More precisely, by \eqref{intersect} we have the primary decomposition
\begin{equation*}
	I(G) = \bigcap \big(\x_{C}=\underset{x_i\in C}{\prod} x_i \mid C \text{ is a minimal vertex cover of } G\big).
\end{equation*}
We therefore get a one-to-one correspondence between the set of minimal vertex covers of a simple graph and the set of associated prime ideals of the corresponding edge ideal.
\[
\left\{
\begin{array}{c}
	\text{Minimal vertex} \\
	\text{covers of a} \\
	\text{simple graph } G
\end{array}
\right\}
\quad
\overset{1:1}{\longleftrightarrow}
\quad
\left\{
\begin{array}{c}
	\text{Associated prime} \\
	\text{ideals of } I(G)
\end{array}
\right\}
\]

Some advantages of this one-to-one correspondence are that we can use it to determine the number of the set of associated prime ideals and easily check if a given prime ideal is an associated prime ideal of the corresponding edge ideal. The main purpose of this section is to determine the number of minimal vertex covers of a perfect binary tree. We have the following result.

\begin{thm}\label{numbercover} Let $T$ be a  perfect binary tree of height $h$, and $m_h$ be the number of minimal vertex covers of $T$. Then, there is a recursive relation
	$$\begin{cases}
	m_0=1,m_1=2,m_2=4,\\
	m_{h+1}=2m_hm_{h-2}^4+m_{h-1}^4-m_{h-2}^8 \text{ for all } h\ge 2.
\end{cases}$$
\end{thm}
\begin{proof} Set $x_1$ the root vertex of $T$ and $x_2, x_3$ are two children of $x_1$. Let $C$ be any minimal vertex cover of $T.$
	
	First, we prove two claims, which hold in any subtree $T_{u}$ with root $x_u$ and height $t$ of $T$, as follows:
	
	{\it Claim 1:} Let $C_u$ be an any the minimal vertex cover of $T_u$. If $x_u\notin C_u$ and $t\geq 2$, then number of minimal vertex covers of $T_u$ is $m_{t-2}^4$.
	
	Indeed, because $x_u\notin C_u$, it follows two children $x_{u_1}, x_{u_2}$ of $x_u$ must be contained in $C_u$ to cover edges $\{x_u, x_{u_1}\}, \{x_u, x_{u_2}\}$. Then, there are four perfect binary subtrees of height $t-2$ whose roots are children of $x_{u_1}, x_{u_2}$. Each perfect binary subtree has $m_{t-2}$ elements in the set of the minimal vertex covers. Thus, the number of minimal vertex covers of $T_u$ is $m_{t-2}^4$.
	
	{\it Claim 2:} If $x_u\in C_u$, $t\geq 1$ and ``father" $x_v$ of $x_u$ does not belong to $C_u$, then number of minimal vertex covers of $T_u$ is $m_{t-1}^2$.
	
	Indeed, because $x_u\in C_u$, it follows $x_{u_1}, x_{u_2}$ are not contained in $C_u$.  However, each $x_{u_i}$ for $i=1,2$ is also a root of a perfect binary subtree of $T_u$ and its height is $t-1$. Thus, the number of minimal vertex covers of each subtree of $T_u$ is $m_{h-1}$. Hence, the number of minimal vertex covers of  $T_u$ is $m_{t-1}^2$, as claimed.
	
	Next, we return to prove the theorem by induction. It is not difficult to check  when $h=0, h=1, h=2$, we get $m_0=1,m_1=2,m_2=4.$
	
	Assume that for any perfect binary tree $T$ with height $h>2$ and $m_h$ we have
	$$m_{h}=2m_{h-1}m_{h-3}^4+m_{h-2}^4-m_{h-3}^8.$$
	
	Let $T$ be a perfect binary tree with height $h+1$ and $C$ be any minimal vertex cover of $T.$ Let  $x_1$ be a root vertex of $T$ and $x_2, x_3$ are two children of $x_1$. We also consider the following two cases:
	
	{\bf Case 1:} If $x_1\in C$. Because  $C$ is minimal, there is at least one vertex two vertices in $x_2, x_3$ in $C$. If $x_2\notin C$, then by Claim $1$, the subtree $T_{x_2}$ has the number of minimal vertex covers is $m_{h-2}^4.$ Moreover, $T_{x_3}$ is the perfect binary subtree has height $h-1$. By induction, its number of minimal vertex covers is $m_h$. Thus,  in this case, there are  $m_hm_{h-2}^4$ the minimal vertex covers of $T$. Similar argument, if $x_3\notin C$ there are more $m_hm_{h-2}^4$ the minimal vertex covers of $T$. If neither of $x_2$ nor $ x_3$ is in $C,$ by Claim $1$ again, we can count more $m_{h-2}^8$ minimal vertex covers of $T$. But, this number is computed twice, thus there exists $2m_hm_{h-2}^4-m_{h-2}^8$ in this case.
	
	{\bf Case 2:} If $x_1\notin C$. When both of $x_2$ and $ x_3$ must to belong to $C$ to cover $\{x_1, x_2\}$ and $\{x_1, x_3\}$. Now, applying the Claim $2$ to two independent subtrees with roots $x_2$ and $x_3$, we have $m_{h-1}^4$ the minimal vertex cover sets.
	
	Combine Case $1$ and Case $2$, we have the number of minimal vertex covers of $T$ is 
	$$m_{h+1}=(2m_hm_{h-2}^4-m_{h-2}^8)+m_{h-1}^4,$$
	as required.
\end{proof}
%Although determining an explicit formula for $F_h$ is very difficult, we can still compute $F_h$ from this recursive formula and clearly see its benefit, for example, we compute the first few values of $F_h$.
%\begin{center}
%\begin{tabular}{|c|c|c|c|}
%	\hline
%	$h$ & $|A_h|$ & $|B_h|$ & $F_h$ \\
%	\hline
%	0 &1 &0 &1 \\
%	\hline
%	1 & 1& 1& 2\\
%	\hline
%	2 &3 &1 &4 \\
%	\hline
%	3 & 7& 16& 23\\
%	\hline
%	4 &480 &256 & 736\\
%	\hline
%\end{tabular}
%\end{center}

Note that because of the complexity of the primary decomposition problem, even when $T$ is a perfect binary tree of height $5$, it is not easy to determine the number of associated prime ideals of $I(T)$. However, thanks to the above recursive formula, we can easily compute the number of associated prime ideals of $I(T$), it is $$|\ass(I(T))|=2.736.4^4+23^4-4^8=591137.$$

We establish a specific formula for the vertex cover number of a perfect binary tree.
\begin{prop}\label{matnumber} Let $T$ be a perfect binary tree of height $h>0$. Then 
	\begin{equation*}
		\alpha(T)=\dfrac{2^{h+2}-(3+(-1)^h)}{6}.
	\end{equation*}
\end{prop}

\begin{proof} Denote $\alpha(h)$ is the vertex cover number of a perfect binary tree of height $h$. We will prove the cardinality of a smallest vertex cover set of $G$ is equal to the total number of vertices at levels $1, 3, 5, \ldots, h-3, h-1$ if $h$ is even or at levels $0, 2, 4, \ldots, h-3, h-1$ if $h$ is odd. Note that the number of vertices at level $k$ is $2^k$. This means that we want to show:
	\begin{itemize}
		\item If $h$ is even, $$\alpha(h) = 2(1+4+\dots+4^{(h-2)/2}) = 2 \frac{4^{h/2}-1}{4-1}= \frac{2^{h+1}-2}{3}.$$
		\item If $h$ is odd, $$\alpha(h) = 1+4+\dots+4^{(h-1)/2}= \frac{4^{(h+1)/2}-1}{4-1} = \frac{2^{h+1}-1}{3}.$$
	\end{itemize}
	
	We prove this by strong induction on the height $h$. The cases where $ h=1$ and $h=2$ are trivial. Assume the statement is true for all perfect binary trees of height $k < h$, for some $h > 2$.
	
Let $T$ be a perfect binary tree of height $h$ with root $r$. Let $T_L$ and $T_R$ be the left and right subtrees of $T$ whose roots are two children of $r$. Then $T_L, T_R$ are perfect binary trees of height $h-1$. Let $C$ be a smallest vertex cover of $T$, that is $\alpha(h)=|C|$. We consider two cases for the root $r$.
	
	\textbf{Case $r \in C$}.  If $r \in C$, the edges from $r$ to its children are covered. The rest of $C$ must form the smallest vertex covers for $T_L$ and $T_R$. The cardinality of the cover is $|C_1| = 1 + \alpha(h-1) + \alpha(h-1) = 1 + 2 \cdot \alpha(h-1)$.
	
	\textbf{Case $r \notin C$}. In this case, two children of $r$ must be in $C$. These children cover all edges up to level 2. The rest of $C$ must cover the four subtrees of height $h-2$ rooted at the grandchildren of $r$. The cardinality is $|C_2| = 2 + 4 \cdot \alpha(h-2)$.
	
	It is clear that $\alpha(h) = \min\{|C_1|, |C_2|\} = \min\{1 + 2 \cdot \alpha(h-1), 2 + 4 \cdot \alpha(h-2)\}$. 
	
	Hence, we consider the case $h$ even first. In this case, $h-1$ is odd and $h-2$ is even, by the inductive hypothesis we get $\alpha(h-1) = \frac{2^{h}-1}{3}$ and $\alpha(h-2) = \frac{2^{h-1}-2}{3}$. So,
	\begin{align*}
		|C_1| &= 1 + 2 \cdot\alpha(h-1) = 1 + 2 \left( \frac{2^h-1}{3} \right) = \frac{2^{h+1}+1}{3}, \\
		|C_2| &= 2 + 4 \cdot \alpha(h-2) = 2 + 4 \left( \frac{2^{h-1}-2}{3} \right)= \frac{2^{h+1}-2}{3}.
	\end{align*}
	Since $\frac{2^{h+1}-2}{3} < \frac{2^{h+1}+1}{3}$, we have $\alpha(h) = |C_2| = \frac{2^{h+1}-2}{3}$ as desired. 
	
	Now consider the case  $h$ odd. 	Then $h-1$ is even and $h-2$ is odd. By the inductive hypothesis $\alpha(h-1) = \frac{2^{h}-2}{3}$ and $\alpha(h-2) = \frac{2^{h-1}-1}{3}$. So,
	\begin{align*}
		|C_1| &= 1 + 2 \cdot \alpha(h-1) = 1 + 2 \left( \frac{2^h-2}{3} \right) = \frac{2^{h+1}-1}{3}, \\
		|C_2| &= 2 + 4 \cdot S(h-2) = 2 + 4 \left( \frac{2^{h-1}-1}{3} \right)= \frac{2^{h+1}+2}{3}.
	\end{align*}
	Since $\frac{2^{h+1}-1}{3} < \frac{2^{h+1}+2}{3}$, we have $\alpha(h) = |C_1| = \frac{2^{h+1}-1}{3}$ as desired. 
\end{proof}

Note that a tree is a bipartite graph, by Lemma \ref{konig}, we have the following result.
\begin{cor}
Let $T$ be a perfect binary tree with height $h>0$.  Then the matching number of $T$ is
\begin{equation*}
\dfrac{2^{h+2}-(3+(-1)^h)}{6}.
\end{equation*}
\end{cor}
From  Lemma \ref{nodes},  Proposition \ref{matnumber}, and the fact $\alpha(T)=n-\beta(T)$, we have the corollary.
\begin{cor}\label{dimvsbeta}
	Let $T$  be a perfect binary tree of height $h>0$ and $I$ is its edge ideal in $S$. Then 
	$$\dim(S/I)=\beta(T)=\dfrac{2^{h+3}-3+(-1)^{h}}{6}.$$
\end{cor}
\section{Projective dimensions and last Betti numbers}
In this section, we investigate the projective dimension and last Betti number of the edge ideal of a perfect binary tree. Let $T$ be a perfect binary tree of height $h$ with the vertex set $V(T)=\{x_1,\ldots,x_n\}$. Let $I\subseteq S=\KK[x_1,\ldots, x_n]$ be the edge ideal of $T$. We know by Lemma \ref{nodes} that $T$ has $2^{h+1}-2$ edges which implies $I$ has a minimal set of generators of $2^{h+1}-2$ elements. This yields that the first total Betti number $\beta_1(S/I)=2^{h+1}-2$. The information on other Betti numbers is generally limited. Based on splitting edges, in \cite{HT07}, Ha and Van Tuyl remarked that if a simple graph $T$ has a vertex of degree 1, say $x$, then the edge formed by $x$ and its neighbor $y$ is a splitting edge of $I$. Hence the Betti numbers of $I$ can be calculated through the Betti numbers of $J=(xy)$, $K=I(T\setminus{xy})$ and $J\cap K$. However, the problem of explicitly computing Betti numbers is still very difficult, even for perfect binary trees. Let $T$ be a tree, we recall a known result in \cite{Ki16, MTV} on the depth of $S/I(T)$ through the maximal independent sets of $T$ as follows.
\begin{lem}[\cite{Ki16, MTV}]\label{depthvsInd}
	Let $T$ be a tree. Then $\depth S/I(T) = q(T)$, where $q(T)$ is the minimum size of a maximal independent set of $T$.
\end{lem}
It is useful to apply this to calculate the depth of the edge ideal of a binary tree.
\begin{prop}\label{depthtree}
	Let $T$ be a perfect binary tree of height $h>0$ and $I$ be the edge ideal of $T$. Then the following assertions holds true.
	\begin{enumerate}
		\item[(1)] If $h=3m$ then $\depth(S/I)=\dfrac{4.8^{m}+3}{7}$.
		\item[(2)] If $h=3m+1$ then $\depth(S/I)=\dfrac{8^{m+1}-1}{7}$.
		\item[(3)] If $h=3m+2$ then $\depth(S/I)=\dfrac{2(8^{m+1}-1)}{7}$.
	\end{enumerate}
\end{prop}

\begin{proof} 
	By Lemma \ref{depthvsInd}, it is sufficient to find the minimum size of a maximal independent set of $T$.
	Let $M$ is the maximal independent set of $T$ with the fewest elements, that is, $|M|$ is the minimum size of a maximal independent set of $T$.\\ 
	$\bullet$ Considering the case $h=3m$, we will show that if $h=3m$ then $M$ consists of the root and all vertices at heights $h-1, h-4, h-7,\ldots$, and $2$. Hence, since there $2^{i}$ vertices at the $i$th height, we have that $$|M|=\sum_{i=1}^{m}2^{3m+2-3i}+1=\dfrac{4.8^{m}+3}{7}.$$
	
	By the definition of $M$, it is clear that no two vertices in the set $M$ are adjacent, which implies that $M$ is an independent set. Moreover, every vertex at heights $h$ and $h-2$ is adjacent to a vertex at height $h-1$, every vertex at heights $h-3$ and $h-5$ is adjacent to a vertex at height $h-4, \ldots $ Finally,  the vertices at heights $3$ and $1$ are adjacent to a vertex at height $2$. Therefore, every vertex not in $M$ must be adjacent to a vertex in $M$, which implies that $M$ is a maximal independent set. 
	
	To see the minimality of $M$, we use induction on $m$. The case  $m=1$, then $|M|=5$. Moreover, a direct computation shows that any maximal independent set must have size at least $5$.
	
	Suppose $m>1$ and the assertion holds for height $3(m-1)$. Let $T$
	be perfect binary of height $h=3m$ and let $M'$ be any maximal independent set of $T$.

	Focus on the bottom two heights $h-1,h$. Each parent at height $h-1$ together with its two children at height $h$ requires at least one chosen vertex from \{that parent, its two children\} in order to dominate the leaves.  
	Since there are $2^{\,h-1}$ such parents, we obtain
	\[
	|M'\cap\{\text{vertices at heights }h-1,h\}| \;\ge\; 2^{\,h-1}.
	\]
	
	Now remove the bottom three heights $h-2,h-1,h$. What remains is a perfect binary tree of height $h-3=3(m-1)$.  
	The restriction of $M'$ to this subtree is a maximal independent set there, so by the induction hypothesis, it has size at least
	$1+\sum_{i=0}^{m-2} 2^{h-4-3i}$ elements.
	
	Combining these bounds yields
	$$ |M'|\ge2^{h-1}+\Big(1+\sum_{i=0}^{m-2} 2^{3m-4-3i}\Big)
	= 1+\sum_{i=0}^{m-1}2^{3m-1-3i}=\frac{4.8^{m}+3}{7}.$$
	Thus, the set consists of the root and all vertices at heights $h-1, h-4, h-7,...,$ and $2$ is the maximal independent set of minimum size. Hence, $$\depth(S/I)=|M|=\frac{4.8^{m}+3}{7}.$$
	$\bullet$ A completely similar argument shows that if $h=3m+1$, the set consists of all vertices at heights $h-1, h-4, h-7,...,$ and $0$ is the maximal independent set of minimum size of $T$. Hence, $$\depth(S/I)=\sum_{i=1}^{m+1}2^{3m+3-3i}=\frac{8^{m+1}-1}{7}.$$ 
	Finally, if $h=3m+2$, the set consists of all vertices at heights $h-1, h-4, h-7,...$ and $1$ is the maximal independent set of minimum size of $T$. Hence, $$\depth(S/I)=\sum_{i=1}^{m+1}2^{3m+4-3i}=\dfrac{2(8^{m+1}-1)}{7}.$$ 
\end{proof}

As an immediate consequence, we determine the projective dimension of edge ideals of perfect binary trees.

\begin{thm}\label{pdim}
	Let $T$ be perfect binary tree of height $h>0$ and $I$ be the edge ideal of $T$. Then the following assertions holds true.
	\begin{enumerate}
		\item[(1)] If $h=3m$ then $\pd_S(S/I)=\dfrac{10(8^m-1)}{7}$.
		\item[(2)] If $h=3m+1$ then $\pd_S(S/I)=\dfrac{20.8^m-6}{7}$.
		\item[(3)] If $h=3m+2$ then $\pd_S(S/I)=\dfrac{40.8^m-5}{7}$.
	\end{enumerate}
\end{thm}
\begin{proof}
	We have by the Auslander-Buchsbaum's formula that 
	$$\pd_S(S/I)=\depth(S)-\depth(S/I)=2^{h+1}-1-\depth(S/I).$$
	Hence, by Proposition \ref{depthtree} we deduce the desired result.
\end{proof}

The following result implies that the total Betti number of the edge ideal of a perfect binary tree at the highest homological degree always equals one.
\begin{thm}\label{main2}
	Let $T$ be a perfect binary tree of height $h>0$, and $I$ its edge ideal in polynomial ring $S=\mathbb K[x_v\mid x_v\in V(T)]$. Denote by $p_h=\pd_S(S/I)$ the projective dimension of $S/I$. Then the last total Betti number of $S/I$ is $\beta_{p_h}(S/I) = 1.$
\end{thm}

\begin{proof} We prove $\beta_{p_h}(S/I) = 1$ by induction on $h$.  If $h=1$, the tree $T$ is the complete bipartite graph $K_{1,2}$ with edge ideal $I = (x_1x_2, x_1x_3)$. The minimal free resolution of $S/I$ is
	$$ 	0 \to S\xrightarrow{\binom{-x_3}{x_2}} S^2\xrightarrow{(x_1x_2,x_1x_3)} S \to S/I \to 0.$$
	So $\beta_1(S/I)=1$ and the assertion holds.
	
	We consider $h\ge 2$ and assume that the statement holds for perfect binary trees of height less than $h$.
	Let $x_1$ be the root of $T$. Let $T'$ and $T''$ be the left and right subtrees obtained from $T$ by removing the root $x_1$ and two incidents with $x_1$. Let $x_{2}$ and $x_{3}$ be the roots of $T'$ and $T''$, respectively. We first show the following claims.
	
	{\it Claim 1. $\pd\left(\dfrac{S}{I+(x_1)}\right)=1+2p_{h-1}$}.
	
	Indeed, since $V(T')\cap V(T'')=\emptyset$, we have a disjoint union $V(T)=V(T')\sqcup V(T'')\sqcup \{x_1\}$. Hence, there is an isomorphism
	\begin{equation}\label{eq9}
		\frac{S}{I+(x_1)}\cong \frac{\KK[x_1]}{(x_1)}\otimes_{\mathbb K}\frac{S'}{I(T')}\otimes_{\mathbb K} \frac{S''}{I(T'')},
	\end{equation}
	where $S'=\KK[x_v\mid x_v\in V(T')]$ and $S''=\KK[x_v\mid x_v\in V(T'')]$. Moreover, since $T'$ and $T''$ are  also perfect binary trees of height $h-1$, we have by the assumption that $\pd_{S'} S'/I(G')=\pd_{S''} S''/I(G'')=p_{h-1}$. Hence, since $\pd\frac{\KK[x_1]}{(x_1)}=1$, we get
	$$
	\pd_S\left(\frac{S}{I+(x_1)}\right)=1+p_{h-1}+p_{h-1}=1+2p_{h-1}.
	$$
	as desired. 
	
	{\it Claim 2. $\pd_S\left(\dfrac{S}{I:x_1}\right)=2+4p_{h-2}$}.
	
	Indeed, we observe that all monomials in minimal set of generators of $I$ that involve $x_1$ are $x_1x_{2}$ and $x_1x_{3}$. So, there is an equality
	$$I:x_1=(x_{2},x_{3})+I(T')+I(T'').$$
	Note that in the quotient ring $\dfrac{S}{I:x_1}$ then $x_{2}=x_{3}=0$. By setting $x_{2}=0$ then $T'$ is splitted into two perfect binary subtrees $T^1, T^2$ of same height $h-2$. Similarly, $T'$ breaks into two perfect binary subtrees $T^3, T^4$ of same height $h-2$ by putting $x_{3}=0$. Moreover, we have a disjoint union $V(T)=V(T^1)\sqcup V(T^2)\sqcup V(T^3)\sqcup V(T^4)\sqcup \{x_1,x_{2},x_{2}\}$ which implies an isomorphism
	\begin{equation}\label{eq10}
		\frac{S}{I:x_1}\cong \frac{\KK[x_1,x_{2},x_{3}]}{(x_{2},x_{3})}\otimes_{\KK}\frac{S^1}{I(T^1)}\otimes_{\KK} \frac{S^2}{I(T^2)}\otimes_{\KK} \frac{S^3}{I(T^3)}\otimes_{\KK} \frac{S^4}{I(T^4)},
	\end{equation}
	where $S^i=\KK[x_v\mid x_v\in V(T^i)]$, $1\le i\le 4.$ Also by the assumption we have $\pd_{S^i} S^i/I(T^i)=p_{h-2}$ for all $1\le i\le 4.$ Hence, since $\pd\frac{\KK[x_1,x_{2},x_{3}]}{(x_{2},x_{3})}=2$, we get
	$$
	\pd_S\left(\frac{S}{I:x_1}\right)=2+4p_{h-2}.
	$$
	
	Now look at the short exact sequence 
	$$0\to \frac{S}{I:x_1}(-1)\xrightarrow{\cdot x_1}\frac{S}{I}\to \frac{S}{I+(x_1)}\to 0.$$
	By applying $\Tor_i^S(-,\mathbb K)$ to this exact sequence, we get the long exact sequence
	\begin{align*}
		\cdots\to\Tor_{i}^S\left(\frac{S}{I:x_1}(-1),\mathbb K\right)\to &\Tor_i^S\left(\frac{S}{I},\mathbb K \right)\to \Tor_i^S\left(\frac{S}{I+(x_1)},\mathbb K\right)\to\\ &\Tor_{i-1}^S\left(\frac{S}{I:x_1}(-1),\mathbb K\right)\to \Tor_{i-1}^S\left(\frac{S}{I},\mathbb K\right)\to\cdots
	\end{align*}
	In particular, we have
	\begin{align*}
		\cdots\to\Tor_{p_h+1}^S\left(\frac{S}{I+(x_1)},\mathbb K\right)\to \Tor_{p_h}^S\left(\frac{S}{I:x_1}(-1),\mathbb K\right)&\to \Tor_{p_h}^S\left(\frac{S}{I},\mathbb K\right)\to\\ \Tor_{p_h}^S\left(\frac{S}{I+(x_1)},\mathbb K\right)&\to \Tor_{p_h-1}^S\left(\frac{S}{I:x_1}(-1),\mathbb K\right)\to\cdots
	\end{align*}
	We consider the following three cases:\\
	$\bullet$ Case $h=3m+1$. Then, by Corollary \ref{pdim} we have $$\begin{cases}
		p_h&=p_{3m+1}=\dfrac{20.8^m-6}{7},\\
		1+2p_{h-1}&=1+2p_{3m}=1+2\dfrac{10(8^m-1)}{7}=\dfrac{20.8^m-13}{7},\\
		2+4p_{h-2}&=2+4p_{3(m-1)+2}=2+4\dfrac{40.8^{m-1}-5}{7}=\dfrac{20.8^m-6}{7}.
	\end{cases}$$
	Hence,
	$p_h=2+4p_{h-2}>1+2p_{h-1}$. Therefore, 
	$$\Tor_{p_h}^S\left(\frac{S}{I+(x_1)},\KK\right)=\Tor_{p_h+1}^S\left(\frac{S}{I+(x_1)},\KK\right)=(0),$$
which yields
	$$\Tor_{p_h}^S\left(\frac{S}{I},\KK\right)\cong \Tor_{p_h}^S\left(\frac{S}{I:x_1},\KK\right).$$
	Moreover, by (\ref{eq10}) and the fact that   $\pd_S\left(\dfrac{S}{I:x_1}\right)=2+4p_{h-2}=p_h$, there are equalities
	\begin{align*}
		\beta_{p_h}(S/I)&=\dim_{\KK}\Tor_{p_h}^S\left(S/I,\KK\right)=\dim_{\mathbb K}\Tor_{2+4p_{h-2}}^S\left(\frac{S}{I:x_1},\KK\right)\\
		&=\beta_2(\frac{\KK[x_1,x_{2},x_{2}]}{(x_{2},x_{2})}).\prod_{i=1}^{4}\dim_{\mathbb K}\Tor_{p_{h-2}}^{S^i}\left(\frac{S^i}{I(T^i)},\KK\right)=1.1.1.1.1=1.
	\end{align*}
	Here the last equality is follows from the inductive hypothesis.\\
	$\bullet$ Case $h=3m+2$. Then, also by Corollary \ref{pdim} we have $$\begin{cases}
		p_h&=p_{3m+2}=\dfrac{40.8^m-5}{7},\\
		1+2p_{h-1}&=1+2p_{3m+1}=1+2\dfrac{20.8^m-6}{7}=\dfrac{40.8^m-5}{7},\\
		2+4p_{h-2}&=2+4p_{3m}=2+4\dfrac{10.8^m-10}{7}=\dfrac{40.8^m-26}{7}.
	\end{cases}$$
	Hence,
	$p_h=1+2p_{h-1}>3+4p_{h-2}$. Therefore, 
	$$\Tor_{p_h-1}^S\left(\frac{S}{I:x_1},\KK\right)=\Tor_{p_h}^S\left(\frac{S}{I:x_1},\KK\right)=(0),$$
	which implies that 
	$$\Tor_{p_h}^S\left(\frac{S}{I},\KK\right)\cong \Tor_{p_h}^S\left(\frac{S}{I+(x_1)},\KK\right).$$
	
	Moreover, by (\ref{eq9}) and the fact that   $\pd_S\left(\dfrac{S}{I+(x_1)}\right)=1+2p_{h-1}=p_h$, there are equalities
	\begin{align*}
		\beta_{p_h}(S/I)&=\dim_{\KK}\Tor_{p_h}^S\left(S/I,\KK\right)=\dim_{\KK}\Tor_{1+2p_{h-1}}^S\left(\frac{S}{I:x_1},\KK\right)\\
		&=\beta_1(\frac{\KK[x_1]}{(x_1)}).\dim_{\KK}\Tor_{p_{h-1}}^{S'}\left(\frac{S'}{I(T')},\KK\right).\dim_{\KK}\Tor_{p_{h-1}}^{S''}\left(\frac{S''}{I(T'')},\KK\right)=1.1.1=1.
	\end{align*}
	Here the last equality is also due to the inductive hypothesis.

	$\bullet$ Finally, we consider the case $h=3m, m>0$. Let $x_2, x_3$ are the children of $x_1$ and $x_4, x_5$ are the children of $x_2$. Let $T_{x_4}$ and $T_{x_5}$ are perfect binary subtrees of height $h-2$ of $T$ with the roots $x_4$ and $x_5$ respectively. Let $G$ is the graph obtained from $T$ by removing $T_{x_4},  T_{x_5}$ and the vertex $x_2$. We will investigate two ideals $I+(x_2)$ and $I:x_2$.
	
	Firstly, it is clear that $$I+(x_2)=(x_2)+I_1+I_2+J,$$
	where  $I_1, I_2, J$ are the edge ideals of $T_{x_4},  T_{x_5}$ and $G$ respectively. Since $V(T)=V(T_{x_4})\sqcup V(T_{x_5})\sqcup V(G)\sqcup \{x_{2}\}$, we have an isomorphism
	\begin{equation}\label{eq11}
		\frac{S}{I+(x_2)}\cong \frac{\KK[x_{2}]}{(x_{2})}\otimes_{\KK}\frac{S^1}{I_1}\otimes_{\KK} \frac{S^2}{I_2}\otimes_{\KK} \frac{S'}{J},
	\end{equation}
	where $S^1=\KK[x_v\mid \x_v\in V(T_{x_4})], S^2=\KK[x_v\mid \x_v\in V(T_{x_5})]$ and $S'=\KK[x_v\mid \x_v\in V(G)]$. 
	
	We first compute the projective dimension of $S'/J$. Let $H$ is are perfect binary subtree of height $h-1$ of $G$ with the roots $x_3$. Let $I_3$ the edge ideal of $H$ and let $x_6, x_7$ are the children of $x_3$ with $S''=\KK[x_v\mid x_v\in V(H)]$. Observe that $H$ has 4 perfect binary subtrees of height $h-3$ whose roots are children of $x_6, x_7$, we call them $H_i \, (1\le j\le 4)$. We have $$I_3:x_1x_3=(x_6, x_7)+J_1+J_2+J_3+J_4.$$
	From the disjoint union $V(G)=\{x_6,x_7\}\sqcup V(H_1)\sqcup V(H_2)\sqcup V(H_3)\sqcup V(H_4)$ and $H_i$ are perfect binary trees of height $h-3$, by similar reasoning as the previous two cases we also have
	$$\pd_{S'}\frac{S'}{I_3:x_1x_3}=2+4p_{h-3}=2+p_{3(m-1)}=\frac{5.8^m-26}{7}.$$
	
	Now we have a short exact sequence 
	$$0\to \frac{S'}{I_3:x_1x_3}\xrightarrow{\cdot x_2x_3}\frac{S'}{I_3}\to \frac{S'}{J}\to 0 \qquad\qquad(\sharp).$$
	Now since $H$ is a perfect binary tree of height $h-1$, $$\pd_{S'}(S'/I_3)=p_{h-1}=p_{3(m-1)+2}=\dfrac{5.8^m-5}{7}\neq \pd_{S'} (S'/I_3:x_1x_3),$$ we get $\pd_{S'}(S'/J)=\pd_{S'}(S'/I_3)=\dfrac{5.8^m-5}{7}.$
	Then by \cite[Theorem 2.1.1]{B} the minimal free resolution of $S'/J$ is obtained by mapping cone applied in $(\sharp)$, we have
	$$\beta_{ij}(S'/J)=\beta_{ij}(S'/I_3)+\beta_{i-1,j-2}(S'/I_3:x_1x_3).$$
	In particular, since $p_{h-1}=\pd_{S'}(S'/J)=\pd_{S'}(S'/I_3)>\pd_{S'} (S'/I_3:x_1x_3)+1$, we get
	$\beta_{p_{h-1}}(S'/J)=\beta_{p_{h-1}}(S'/I_3)=1$.
	The last equality follows from the inductive hypothesis.

	Now the isomorphism (\ref{eq11}) implies that   	\begin{align*}
		\pd_S\frac{S}{I+(x_2)}&=1+\pd_{S^1}\frac{S^1}{I_1}+\pd_{S^2}\frac{S^2}{I_2}+\pd_{S'}\frac{S'}{J}\\
		&=1+2p_{h-2}+p_{h-1}=1+2p_{3m-2}+p_{3m-1}\\
		&=\frac{10.8^m-10}{7}=p_h.
	\end{align*}
	Hence, again by (\ref{eq11}), we have equalities $$\beta_{p_h}(S/(I+(x_2)))=\beta_1(\KK[x_2]/(x_2)).\beta_{p_{h-2}}(S^1/I_1).\beta_{p_{h-2}}(S^2/I_2).\beta_{p_{h-1}}(S'/J)=1.1.1.1=1.$$
	
	On the other hand, we have $$I:x_2=(x_1,x_4,x_5)+\mathfrak{a}_1+\mathfrak{a}_2+\mathfrak{a}_3+\mathfrak{a}_4+I_3,$$
	where $\mathfrak{a}_i \, (1\le i\le 4)$ be edge ideals of perfect binary subtrees of height $h-3$ whose roots are children of $x_4$ and $x_5$. Hence, by the same arguments as the previuos two
	cases we can show that
	\begin{align*}
		\pd_S(S/I:x_2)&=3+4p_{h-3}+p_{h-1}=3+4p_{3(m-1)}+p_{3(m-1)+2}\\
		&=3+4.\frac{10.8^{m-1}-10}{7}+\frac{40.8^{m-1}-5}{7}\\
		&=\frac{10.8^m-24}{7}.
	\end{align*}
	So, $$\pd_S(S/I)=\pd_S(S/(I+(x_2)))=p_h=\frac{10.8^m-10}{7}>\pd_S(S/(I:(x_2)))+1.$$
	Hence, the long exact sequence induced by the exact sequence
	$$0\to \frac{S}{I:x_2}\xrightarrow{\cdot x_2}\frac{S}{I}\to \frac{S}{I+(x_2)}\to 0,$$
	implies that $$\Tor^S_{p_h}(S/I,\KK)\cong \Tor^S_{p_h}(S/I+(x_2),\KK).$$
	In particular, $\beta_{p_h}(S/I)=\beta_{p_h}(S/I+(x_2))=1$ as desired. We therefore have finished the proof.
\end{proof}

\begin{rem} Note that Corollary \ref{dimvsbeta} and Theorem \ref{main2} imply that if $T$ is a perfect binary tree of height $h > 0$ then $S/I$ is not a Cohen-Macaulay ring. It has been well-known that the last Betti number is the type of $S/I$ that is defined by $\dim_\KK{\rm Ext}^t(\KK, S/I)$ with $t=\depth_SS/I$. Moreover, a Cohen-Macaulay ring is a Gorenstein if and only if it has type 1.  Therefore, Theorem \ref{main2} provides a class of non-Gorenstein rings of type 1. 
\end{rem}

\section*{Acknowledgments}
We are grateful to Prof. Tran Nam Trung for his valuable comments.
The first author is partially supported by Vietnam National Foundation for Science and Technology Development (Grant \#101.04-2024.07).% The third author is supported by the Ministry of Education and Training of Vietnam through the Grant number B.2026-SP2-01.

\vspace{1cm}
\noindent {\bf Data Availability} Data sharing is not applicable to this article as no datasets were generated or analyzed during the current study.

\noindent {\bf Conflict of interest} There are no competing interests of either financial or personal nature.


\begin{thebibliography}{99} 
\bibitem{B}  R. Bouchat, {\it Free resolutions of some edge ideals of simple graphs}, J. Commut. Algebra.,  {\bf 2(1)} (2010), 1-36.


\bibitem{BH} W. Brun and J. Herzog,  Cohen-Macaulay rings, {\it Cambridge Studies in Advanced Mathematics} 39, Cambridge University Press, Cambridge 1993.

\bibitem{BM} J. A. Bondy, U. S. R. Murty, Graph Theory, {\it Graduate Texts in Mathematics},  {\bf 244},  Springer New York, 2008.
\bibitem{Co84} D. Cohen, {\it Counting stable sets in trees}, in Seminaire Lotharingien de combinatoire, R. K\"{o}nig, ed., Institute de Recherche Math\'{a}matique Avanc\'{e}e pub., Strasbourg, France, 1984, 48--52.

\bibitem{D} R. Diestel, {\it Graph theory}, 2nd. edition, Springer: Berlin/Heidelberg/New York/Tokyo, 2000.

%\bibitem[EG]{EG} D. Eisenbud, M. Green, K. Hulek, S. Popescu, {\it Restricting linear syzygies: Algebra and geometry}, Compos. Math. 141 (2005) 1460--1478.

%\bibitem[EK90]{EK90} S. Eliahou, M. Kervaire, {\it Minimal resolutions of some monomial ideals,} J. Algebra, {\bf 129} (1990), 1–25.

%\bibitem[Fat01]{Fat01} G. Fatabbi, {\it On the resolution of ideals of fat points}, J. Algebra, {\bf 242} (2001), 92--108.

%\bibitem[FHT] {FHT} C. A. Francisco, H. T. H{$\rm \grave{a}$}, A. V. Tuyl (2010), {\it Associated primes of monomial ideals and odd holes in graphs}, J Algebr Comb, {\bf 32}, 287--301.


\bibitem{Fu} Z. F\"{u}redi, {\it The number of maximal independent sets in connected graphs}, J. Graph Theory, {\bf 11} (1987), 463--470.

\bibitem{HT07} H. T.  Ha, A. V. Tuyl, {\it Splittable ideals and the resolutions of monomial ideals}, J. Algebra,  {\bf 309} (2007), 405–425.

\bibitem{HSM} H. T.  Ha, S. E. Morey, {\it Embedded associated primes of powers of square-free monomial ideal}, J. Pure Appl. Algebra,  {\bf 214(4)} (2010), 301--308.

%\bibitem[HHT1] {HHT1} J. Herzog, T. Hibi and N. V. Trung, {\it Symbolic powers of monomial ideals and vertex cover algebras},  Adv. Math. {\bf 210}(1) (2007), 304 - 322.

\bibitem{HM} A. Hirano, K. Matsuda,  {\it Matching numbers and dimension edge ideals},  Graphs and Combinatorics, {\bf 37(3)} (2021), 761--774.


\bibitem{HT19} D. T. Hoang and T. N. Trung, {\it Coverings, matchings and the number of maximal independent sets of graphs}, Australas. J. Combin., vol. 73, {\bf3} (2019), 424--431.

\bibitem{HT} M. Hujter and Z. Tuza, {\it The number of maximal independent sets in triangle-free graphs}, SIAM J. Discrete Math., {\bf 6} (1993), 284--288.

\bibitem{Ki16} K. Kimura, {\it Non-vanishingness of Betti numbers of edge ideals and complete bipartite subgraphs}, Commun. Algebra, {\bf 44} (2016), 710--730.


\bibitem{Liu} J. Liu, {\it Maximal independent sets in bipartite graphs}, J. Graph Theory, {\bf 17} (1993), 495--507.

%\bibitem{MS} E. Miller and B. Sturmfels, {\it Combinatorial commutative algebra}, Springer, 2005.

\bibitem{MTV} N. C. Minh, T. N. Trung, and T. Vu, Depth of powers of edge ideals of cycles and trees. https://doi.org/10.48550/arXiv.2308.00874.


\bibitem{Sa88} B. E. Sagan, {\it A note on independent sets in trees}, SIAM J. Discrete Math., {\bf 1} (1988), 105--108.

\bibitem{Wi86} H. S. Wilf, {\it The number of maximal independent sets in a tree}, SIAM J. Algebraic Discrete Methods, {\bf 7} (1986), 125--130.

	
\end{thebibliography}
\end{document}